\documentclass[10pt]{birkjour}
\newtheorem{thm}{Theorem}[section]

\newtheorem{lem}[thm]{Lemma}
\newtheorem{prop}[thm]{Proposition}
\theoremstyle{definition}

\theoremstyle{remark}

\numberwithin{equation}{section}
\begin{document}
	
	%
	%
	%
	%
	%

	%
	\author[Jafar Ahmed,  Ajebbar Omar, Elqorachi Elhoucien]{Jafar Ahmed$^{1}$, Ajebbar Omar$^{2}$, Elqorachi Elhoucien$^{1}$}
	
	\address{%
		$^1$
	 Ibn Zohr University, Faculty of sciences,
		Department of mathematics,\\
		Agadir,
		Morocco}
	\address{%
	$^2$
	Sultan Moulay Slimane University, Multidisciplinary Faculty,
	Department of mathematics and computer science,
	Beni Mellal,
	Morocco}
\email{hamadabenali2@gmail.com, omar-ajb@hotmail.com, elqorachi@hotmail.com }
\thanks{2020 Mathematics Subject Classification: Primary 39B52, Secondary 39B32}
\thanks{Keywords: Functional equations, Semigroups, Multiplicative function, Additive function, Cosine-Sine equation.}

	\title[A Kannappan-cosine functional equation on semigroups]{A Kannappan-cosine functional equation on semigroups}
	\begin{abstract} In this paper
		we determine the complex-valued solutions of the Kannappan cosine functional law $g(xyz_{0})=g(x)g(y)-f(x)f(y)$, $x,y\in S$,  where $S$ is a semigroup and $z_{0}$ is a fixed element in $S.$
	\end{abstract}
\maketitle
\section{Introduction}
The addition law for cosine is $$cos(x+y)=cos(x)\,cos(y)-sin(x)\,sin(y),\,x,y\in\mathbb{R}.$$
It is at the origin of the following functional equation on any semigroup $S$:
\begin{equation}\label{1}
g(xy)=g(x)g(y)-f(x)f(y),\,x,y \in S
\end{equation}
for the unknown functions $f,g:S\to \mathbb{C}$,
which is called  cosine addition law. In Acz\`{e}l's monograph \cite[Section 3.2.3]{a} we find the early results about (\ref{1}) and its continuous real valued solutions for $S=\mathbb{R}$.\\
The functional equation (\ref{1}) has been solved on groups by Poulsen and Stetk\ae r (\cite{p.h.}), on semigroups generated by their squares by Ajebbar and Elqorachi (\cite{a.e}), and recently by Ebanks \cite{b} on   semigroups \\
In \cite[ Theorem 3.1]{h} Stetk\ae r solved the  following functional equation
\begin{equation}\label{6-2}
g(xy)=g(x)g(y)-f(y)f(x)+\alpha f(xy), \hspace{0.4cm}x,y \in S,
\end{equation}
where $\alpha$ is a fixed constant in $\mathbb{C} $. He expressed the solutions in terms of multiplicative functions and the solution of the special case of sine addition law. In  \cite[Proposition 16]{i} he solved the functional equation
\begin{equation}\label{4}
f(xyz_{0})=f(x)f(y),\,x,y \in S,
\end{equation}
on semigroups, and where $z_{0}$ is a fixed element in $ S$. We shall use these results in our computations.\\ In this paper we deal with the following Kannappan-cosine addition Law
\begin{equation}\label{2}
g(xyz_{0})=g(x)g(y)-f(x)f(y),\,x,y \in S
\end{equation}
on a semigroup $S$. The functional equation (\ref{2}) is called Kannappan functional equation because it  brings up   a fixed element $z_{0}$ in $S$ as in the paper of Kannappan \cite{kan}. \\
In the special case where $\{f,g\}$ is linearly dependent and $g\neq0$ we get that there exists a constant $\lambda\in\mathbb{C}$ such that the function $(1-\lambda^{2})g$ satisfies the functional equation (\ref{4}).\\
If $S$ is a monoid with identity element $e$, and $f(e)=0$ and $g(e)\neq 0$, or $g(e)=0$ and $f(e)\neq 0$, the last functional equation is the cosine addition law which was solved  recently on general semigroups by Ebanks \cite{b}.\\Now, if $\alpha:=f(e)\neq0$ and $\beta:g(e)\neq0$ we get that the pair $\left(\dfrac{g}{\beta},\dfrac{f}{\beta}\right)$ satisfies the following functional equation \begin{equation*}\dfrac{g}{\beta}(xy)=\dfrac{g}{\beta}(x)\dfrac{g}{\beta}(y)-\dfrac{f}{\beta}(x)\dfrac{f}{\beta}(y)+\dfrac{\alpha}{\beta}\dfrac{f}{\beta}(xy),\end{equation*}
which is of the form (\ref{6-2}), and then explicit formulas
for $f$ and $g$ on groups exist in the literature (see for example \cite[Corollary 3.2.]{zeg}).\\
The natural general setting of the functional equation (\ref{2}) is for $S$ to be
a semigroup, because the very formulation of (\ref{2}) requires only an associative
composition in $S$, not an identity element and inverses. Thus we study in the
present paper Kannappan-cosine functional equation (\ref{2}) on semigroups $S$, generalizing previous works in which $S$ is a group.
So, the result of the present paper is a naturel continuation of results contained in the literature.
\\\\The purpose of the present paper is to show how  the relations
between (\ref{2}) and (\ref{6-2})-(\ref{4}) on monoids extends to the much wider framework, in which $S$ is a semigroup. We find explicit formulas for the solutions, expressing them in terms of
homomorphisms and additives
maps from  the semigroup $S$ into $\mathbb{C}$ (Theorem 4.1). The continuous solutions on topological semigroups are also found.
\section{Set up, notations and terminology}
Throughout this paper $S$ is a semigroup (a set with an associative composition) and $z_0$ is a fixed element in $S$.
\\If $S$ is topological, we denote by $\mathcal{C}(S)$ the algebra of continuous functions from $S$ to field of complex numbers $\mathbb{C}$.\\
Let $f : S\longrightarrow \mathbb{C}$ be a function. We say that $f$ is central if $f(xy) = f(yx)$ for all $x,y \in S$, and that $f$ is abelian if
$ f(x_1 x_2,...,   x_n) = f(x_{\sigma(1)}x_{\sigma(2)},...,x_{\sigma(n)})$
for all $x_1, x_2, . . . , x_n \in S$, all permutations $\sigma$ of $n$ elements and all $n \in \mathbb{N}$.
A map $A:S\longrightarrow \mathbb{C}$ is said to be additive if $A(xy)=A(x)+A(y),\hspace{0.3cm} \text{for} \hspace{0.1cm} \text{all} \hspace{0.1cm} x,y \in S
$ and a map $\chi :S\longrightarrow \mathbb{C}$ is multiplicative if $\chi(xy)=\chi(x)\chi(y),\hspace{0.3cm} \text{for} \hspace{0.1cm} \text{all} \hspace{0.1cm} x,y \in S.$\\
If $\chi \neq 0,$ then the nullspace $I_{\chi} := \{x\in S\hspace{0.1cm} | \hspace{0.1cm} \chi(x) = 0\}$
is either empty or a proper subset of S and $I_{\chi}$ is a two sided ideal in S if not empty and $S\setminus I_{\chi}$ is a subsemigroup of $S$.
\\ Note that additive and multiplicative functions are abelian.\\
For any  subset $T \subseteq S$ let  $T^2 :=\{xy  \hspace{0.1cm}| \hspace{0.1cm}x,y\in T\}$ and for any fixed element $z_0$ in S we let
 $T^2z_0:=\{xyz_0 \hspace{0.1cm}|\hspace{0.1cm} x,y\in T  \}$.\\
To express solutions of our functional equations studied in this paper we will use the set
$P_{\chi}:=\{p\in I_{\chi}\setminus I_{\chi}^2 \hspace{0.1cm} | \hspace{0.1cm}up,pv,upv \in  I_{\chi}\setminus I_{\chi}^2 \hspace{0.1cm}\text{for}\hspace{0.1cm} \text{all} \hspace{0.1cm} u,v \in S\setminus I_{\chi}\}$.
For more details about $P_{\chi}$ we refer the reader to \cite{ c}, \cite{ b} and \cite{d}.			         	 		
\section{Preliminaries}
In this section we give useful results to solve the functional equation (\ref{4}).
 \begin{lem}
		 	  	 	\label{ui}
		 	  	   Let $S$ be  a semigroup, $n\in \textit{N},$
		 	  	   and $\chi,\chi_1,\chi_2,...\chi_n:S\longrightarrow \mathbb{C}$    be different non-zero multiplicative functions. Then  \\
		 	  	    (a) $\{\chi_1,\chi_2,...\chi_n\}$is linearly independent.  \\
		 	  	   (b) If  $A: S\setminus I_{\chi}\longrightarrow \mathbb{C}$ is  a non-zero additive function, then    the set   $\{\chi A, \chi\}$ is linearly independent on $S\setminus I_{\chi}$.
		 	  	 \end{lem}
\begin{proof}
(a): See \cite [Theorem 3.18]{g}. (b): See \cite[Lemma 4.4]{ajjj}.
\end{proof}
The proposition below gives the solutions of the functional equation
\begin{equation}\label{bb}f(xyz_0)=\chi (z_0)f(x)\chi(y)+ \chi (z_0)f(y)\chi(x), \;x,y \in S,\end{equation}
\begin{prop}\label{t2}
Let $S$ be a semigroup, and   $\chi :S\longmapsto\mathbb{C}$ be a multiplicative function such that $\chi(z_0)\neq0$.   If  $f :S\longmapsto\mathbb{C}$ is a solution of (\ref{bb}), then
\begin{equation}\label{15}
f(x)=
\begin{cases}
\chi(x)(A(x)+A(z_0)) \hspace{0.5cm} for\hspace{0.3cm} x\in S \backslash I_{\chi} \\
\rho(x)\hspace{3.2cm} for\hspace{0.3cm} x\in    P_{\chi} \\
0 \hspace{3.8 cm}for\hspace{0.4cm} x\in I_{\chi}\setminus P_{\chi},
\end{cases}
\end{equation}
where $A: S \setminus$
$I_{\chi}\longrightarrow \mathbb{C}$
is additive and    $\rho: P_{\chi}\longrightarrow \mathbb{C}$  is the restriction of $f$ to $P_{\chi} $.\\In addition,    $f$  is abelian and  satisfies the following conditions:\\
(I) $f(xy) = f(yx)=0 $ for all $ x \in I_{\chi} \setminus P_{\chi} $ and $ y \in S \setminus I_{\chi}$.\\
(II) If $x \in \{up, pv, upv\}
$ with $p \in P_{\chi} $ and $u,v \in S \setminus I_{\chi}$, then  $x \in P_{\chi} $ and we have respectively
$\rho(x)=	\rho(p) \chi(u)$, $\rho(x)=	\rho(p) \chi(v)$ or $\rho(x)=	\rho(p) \chi(uv)$.\\Conversely, the  function f  of the form (\ref{15}) define a solution of (\ref{bb}).\\
Moreover, if $S$ is a topological semigroup and $f \in \mathcal{C}(S)$, then   $\chi\in \mathcal{C}(S), A \in\mathcal{C}(S \setminus  I_{\chi})$
and $\rho \in \mathcal{C}(P_{\chi})$.
\end{prop}
\begin{proof}See \cite[Proposition 4.3]{jae}.
\end{proof}
To shorten the way to the solutions of functional equation  (\ref{4}), we prove the following Lemma that contains some key properties.
\begin{lem}
\label{nnn}
Let $S$ be a semigroup and let $f, g :S\longmapsto\mathbb{C}$ the solutions of the functional equation (\ref{4}) with $g\neq 0$. Then\\
(i)	If $f(z_0)=0 $ then	\\	 		 		
(1)
\begin{equation}
\label{0.1}
g(z_0^2)g(xy)= g(z_0)[g(x)g(y)-f(x)f(y)]+f(z_0^2)f(xy),\;x,y\in S.
\end{equation}		 		 		
(2)
\begin{equation*}
g(z_0^2) ^2= g(z_0)^3+ f(z_0^2)^2.
\end{equation*}		 		 		
(3) If $f$ and $g$ are linearly independent then $ g(z_0)\neq0$.\\
		 		 		\\
(ii)	If $f(z_0)\neq0$, then there exist $\mu \in \mathbb{C}$ such that
\begin{equation}
\label{03}
f(xyz_0)=f(x)g(y)+f(y)g(x)+\mu f(x)f(y),\;x,y\in S.
\end{equation}		 		 	
\end{lem}
\begin{proof}
(i)   Suppose that $f(z_0)=0 $.\\
(1)	Making the substitution $(xy,z_0^2)$ and $(xyz_0,z_0)$ in (\ref{4}) we get respectively
$g(xyz_0^3)   =  g(xy) g(z_0^2)-f(xy)f(z_0^2)$,
and
$g(xyz_0^3)=g(xyz_0)g(z_0)-f(xyz_0)f(z_0)=g(z_0)g(x)g(y)-g(z_0)f(x)f(y).$
Comparing these expressions,  we deduce that
		 		 	$
		 		 		g(xy)g(z_0^2)-f(z_0^2)f(xy)=g(z_0)g(x)g(y)-g(z_0)f(y)f(x)
		 		 		$.  This proves  the  identity desired.\\
		 		 		(2) Follows directly by putting $x=y=z_0$ in the equation (\ref{0.1}).\\
		 		 			(3) By contradiction we suppose that $g(z_0)=0$. Then by using (\ref{4}) we get $g(xyz_0^2)= g(x)g(yz_0)-f(x)f(yz_0)=g(xy)g(z_0)-f(xy)f(z_0)=0$
		 		 		since $f(z_0)=g(z_0)=0$. Then we deduce that
		 		 		\begin{equation}
		 		 		\label{fd}
		 		 		g(x)g(yz_0)=f(x)f(yz_0), \; x,y \in S.
		 		 		\end{equation}
		 		 		
		 		 		If $g(yz_0)=0 $ for all $y\in S$ then
		 		 		$
		 		 		0= g(xyz_0)=g(x)g(y)- f(x)f(y), \; x,y \in S.
		 		 		$
		 		 		So,
		 		 		$
		 		 	 g(x)g(y)=f(x)f(y), \;x,y \in S.
		 		 		$
		 		 		Hence  $f=g$ or $f=-g$, which contradicts the fact that f and g are linearly independent.
		 		 	 So  $g\neq0$ on $ S z_0$ then from (\ref{fd}) we get that $g=c_1f$  with $c_1=:f(az_0)/g(az_0)$  for some $a \in S$ such   that  $g(az_0)\neq 0$. This is  also a contradiction, since   $f$ and g are linearly independent.
		 		 		 So we conclude that $g( z_0)\neq0
		 		 		 $. \\(ii) Suppose that $f(z_0)\neq0$.	Making the substitution   $(x,yz_0^2)$ and $(xyz_0,z_0)$ in (\ref{4}) we get respectively
		 		 		$
		 		 			g(xyz_0^3) =  g(x) g(yz_0^2)-f(x) f(yz_0^2)= g(z_0)g(x)g(y)-g(x)f(z_0)f(y) -f(x)f(yz_0^2)
		 		 	$
		 		 		and
		 		 		$
		 		 			g(xyz_0^3) =g(xyz_0 )g(z_0)-f(xyz_0)f(z_0)=g(z_0)g(x)g(y)-g(z_0)f(x)f(y) -f(xyz_0)f(z_0).
		 		 		$
		 		 		Then, by  the associativity of the operation in $S$ we get
		 		 		\begin{equation}
		 		 		\label{00}
		 		 		f(z_0)[f(xyz_0)-f(x)g(y)-f(y)g(x)]=f(x)[f(yz_0^2)-f(y)g(z_0)-f(z_0)g(y)].
		 		 		\end{equation}

		 		 		Since $f(z_0)\neq0$, dividing (\ref{00}) by $f(z_0)$ we get that
		 		 		\begin{equation}
		 		 		\label{999}
		 		 		f(xyz_0)= f(x)g(y)+f(y)g(x)+f(x)\psi(y),
		 		 		\end{equation}
		 		 		where
		 		 	$
		 		 		\psi(y):=f( z_0)^{-1}[f(yz_0^2)-f(y)g(z_0)-f(z_0)g(y)].$
		 		 		Substituting (\ref{999})  come back into (\ref{00}), we find
		 		 		$
		 		 		f(z_0)[f(x)\psi(y)]=f(x) [f(y)\psi(z_0) ],
		 		 	$
		 		 		which implies that
		 		 		$
		 		 		\psi(y)  =\mu f(y)$ with $ \mu:= \dfrac{\psi(z_0)}{f(z_0) }.$
		 		 		Therefore (\ref{999}) becomes
		 		 		$
		 		 		f(xyz_0)= f(x)g(y)+f(y)g(x)+\mu f(x)f(y)
		 		 	$.
This finishes the proof of Lemma \ref{nnn}.
		 		 	\end{proof}
		 		 	\vspace{0.2cm}
\section{Main results}		 		 	
\par Now, we are ready to describe the solutions of the functional equation (\ref{4}).\\
Let $\Psi_{A\chi,\rho}:S\longrightarrow\mathbb{C}$ denote a function $f$ of the form of \cite[Theorem 3.1 (B)]{d}, where $\chi:S\longrightarrow
		 		         	\mathbb{C}$ is  a non-zero multiplicative function, $A:S\setminus I_{\chi}\longrightarrow
		 		         	\mathbb{C}$ is additive,  $\rho:P_{\chi}\longrightarrow
		 		         	\mathbb{C}$ is the restriction of $f$,   and conditions (i)  and (iI) hold.
\begin{thm}
		 		 		\label{jh}
		 		 		The solutions $f,g:S\longrightarrow \mathbb{C}$ of  the functional equation (\ref{4}) are the following  pairs of  functions.\\
		 		 		(1) $f=g=0$.\\
		 		 		(2)   $S\neq S^2z_0$     and we have
		 		 		\begin{equation*}
		 		 		f=\pm g \;  and  \;
		 		 		g(x)=
		 		 		\begin{cases}
		 		 		g_{z_0}(x)\; for\; x\in S  \setminus S^2z_0\\
		 		 		0\hspace{1.5cm} for\hspace{0.3cm} x\in   S^2z_0,
		 		 		\end{cases}
		 		 		\end{equation*}
		 		 	where   $g_{z_0}:S \setminus S^2z_0 \longrightarrow \mathbb{C}$  is an arbitrary non-zero function.\\
		 		 		(3) There exist a constant $d\in \mathbb{C} \setminus\{\pm1\}$ and a multiplicative function
		 		 		$\chi$ on $S$ with $\chi (z_0) \neq0$, such that
		 		 	
		 		 		$f  =\dfrac{d \,\chi(z_0)}{1-d^2}\chi$ and	$ g=\dfrac{\chi(z_0)}{1-d^2}\chi.$
		 		 		\\
		 		 		(4) There exist a constant $c\in \mathbb{C}^{*} \setminus\{ \pm i\}$ and  two  different multiplicative functions
		 		 		$\chi_1$ and $\chi_2$ on $S$, with $ \chi_1(z_0) \neq0 $,  $ \chi_1(z_0)   \neq0$ and $\chi_1(z_0)\chi_1  \neq \chi_2(z_0)\chi_2$
		 		 		such that
		 		 	$$
		 		 		f =-\dfrac{ \chi_1(z_0)\chi_1- \chi_2(z_0)\chi_2 }{ i(c^{-1}+c) }\,\,\text{and}\,\,g =\dfrac{ c^{-1}\chi_1(z_0)\chi_1+c\chi_2(z_0)\chi_2 }{ c^{-1}+c }.  	
		 		 		$$
		 		 			(5)  There exist  constants $q,\gamma\in \mathbb{C}^{*}$ and  two  different non-zero multiplicative functions $\chi_1$ and $\chi_2$ on $S$, with 	$\chi_1(z_0)= \dfrac{q^2-(1+\xi)^2}{2\gamma q} $, $\chi_2(z_0)= -\dfrac{q^2-(1-\xi)^2}{2\gamma q} $ and $\xi:=\pm\sqrt{1+q^2 }$
		 		 			such that
		 		 			$$
		 		 		  f = \dfrac{\chi_1 +\chi_2}{2\gamma}+  \xi
		 		 		  \dfrac{\chi_1 - \chi_2}{2\gamma }	\,\,\text{and}\,\, g=q \dfrac{\chi_1 -\chi_2}{2\gamma}.
		 		 			$$
		 		 		(6) There exist  constants $q\in \mathbb{C}\setminus\{\pm \alpha\}$, $\gamma\in \mathbb{C}^{*}\setminus\{\pm \alpha\}$ and $\delta\in \mathbb{C}\setminus\{\pm 1\}$, and two  different non-zero multiplicative
		 		 		functions
		 		 		$\chi_1$ and $\chi_2$ on $S$, with 	$\chi_1(z_0)=   \dfrac{(1+\delta)^2-(\alpha+q)^2}{2\gamma(1+\delta)} $, $  \chi_2(z_0)=  \dfrac{(1-\delta)^2-(\alpha-q)^2}{2\gamma(1-\delta)} $  and $\delta :=\pm\sqrt{1+q^2-\alpha^2} $
		 		 		such that
		 		 		$$
		 		 		f=\alpha \dfrac{\chi_1 +\chi_2}{2\gamma}+q\dfrac{\chi_1 - \chi_2}{2\gamma} \,\,\text{and}\,\, g = \dfrac{\chi_1 +\chi_2}{2\gamma}+ (1+\delta )
		 		 		\dfrac{\chi_1 - \chi_2}{2\gamma}.
		 		 	$$
		 		 		(7) There exist a constant  $  \beta \in \mathbb{C}^{*} $,  a non-zero multiplicative  function $\chi$
		 		 		on $S$, an additive function $A: S \setminus  I_{\chi}\longrightarrow \mathbb{C}$ and a function $\rho: P_{\chi}\longrightarrow \mathbb{C}$ with	
		 		 			$\chi(z_0)=\dfrac{1  }{\beta } $   and  $A(z_0)=  0$ such that
		 		 		$$
		 		 		f=\dfrac{1}{\beta} \Psi_{A\chi,\rho}\,\,\text{and}\,\,g=\dfrac{1}{\beta} \left( \chi \pm \Psi_{A\chi,\rho}\right).$$
		 		 	 (8) There exist a multiplicative  function $\chi$ on $S$  with $\chi(z_0)\neq0$, an
		 		 		additive function $A:$ $S \setminus  I_{\chi}\longrightarrow \mathbb{C}$ and a function $\rho: P_{\chi}\longrightarrow \mathbb{C}$
		 		 		such that
		 		 		$$
		 		 		f=    A(z_0) \chi   +\Psi_{A\chi,\rho}\,\,\text{and}\,\, g=\left(\chi(z_0) \pm A(z_0)\right) \chi   +\Psi_{A\chi,\rho}.
		 		 		$$ Moreover if $S$ is a topological semigroup and $f \in \mathcal{C}(S)$  then   $g      \in \mathcal{C}(S)$ in cases (1), (2) and (4)-(8), and in (3) if $d\neq0$.
\end{thm}
		 		
\begin{proof}
		 		 		If $g=0$, then (\ref{4}) reduces to $f(x)f(y)=0$ for all $x,y \in S$. This implies that $f=0$, so we get the first part of solutions. From now we may assume that $g\neq0$. \\If $f$ and $g$ are linearly dependent, then there exist $d \in \mathbb{C}$ such that $f=dg$. Substituting this into (\ref{4}) we get the following equation
		 		 		$g(xyz_0)=(1-d^2)g(x)g(y)  \hspace{0.3cm} x,y \in S$. If $d^2=1$, then $g(xyz_0)=0$  for all $x,y \in S$. Therefore $S\neq S^2 z_0$ because $g\neq0$. So we are in solution family (2) with $g_{z_{0}}$ is an arbitrary non-zero function. If $d^2\neq 1$, then by \cite[Proposition 16]{i} there exists a multiplicative function $\chi$ on $S$ such hat $\chi(z_0)\chi:=(1-d^2)g$ and $\chi(z_0)\neq 0$. Then we deduce that
	$ g=\dfrac{\chi(z_0)}{1-d^2}\chi$  and  $f =d\,g=\dfrac{d\,\chi(z_0) }{1-d^2} \chi$, so we have the solution family (3).\\ For the rest of the proof we assume that $f$ and $g$ are linearly independent.
	We split the proof into two cases according to whether
		 		 		$f(z_0)=0$ or $f(z_0)\neq0$.\\
		 		 		{Case 1}. Suppose $f(z_0)=0$. Then by Lemma \ref{nnn} (i)-(3) and (i)-(1)  we have respectively $g(z_0)\neq0$ and
		 		 		\begin{equation}
		 		 		\label{000}
		 		 		g(z_0^2)g(xy)= g(z_0) g(x)g(y)-g(z_0)f(x)f(y) +f(z_0^2)f(xy),\;x,y\in S.
		 		 		\end{equation} {Subcase 1.1}. Assume that  $g(z_0^2)=0$. Then by Lemma \ref{nnn} (i)-(2) and  (i)-(3) we get $ f(z_0^2)\neq0$  since $f$ and $g$ are linearly independent and  then (\ref{000}) can be rewritten as $f(xy)= \gamma f(x)f(y)-\gamma g(x)g(y),$ $x,y \in S,$	where $\gamma:=\dfrac{g( z_0)}{f(z_0^2)}\neq0
		 		 		$. Consequently the pair $(\gamma f, \gamma g)$ satisfies the cosine addition formula (\ref{2}). So, according to   \cite[Theorem 6.1]{h} and  taking into account that  $f$ and $g$ are linearly independent, we know that there are only the  following possibilities.
		 		 		\par (1.1.i)  There exist a constant $q\in \mathbb{C}^{\ast}$ and two  different non-zero multiplicative functions $\chi_1$ and $\chi_2$ on $S$ such that
		 		 		$
		 		 		\gamma g= q\dfrac{\chi_1-\chi_2}{2 } $ and $ \gamma f=\dfrac{\chi_1+\chi_2}{2 }\pm\left(   \sqrt{1+q^2 }\right) \dfrac{\chi_1-\chi_2}{2 },
		 		 		$
		 		 		which gives
		 		 		$
		 		 	  f=\dfrac{\chi_1+\chi_2}{2\gamma }\pm\left( \sqrt{1+q^2 }\right) \dfrac{\chi_1-\chi_2}{2\gamma }$ and $g= q\dfrac{\chi_1-\chi_2}{2 \gamma} $. By putting  $\xi:=\pm\sqrt{1+q^2 }$
		 		 	 and using (\ref{4}) we get
		 		 	$$\dfrac{1}{4\gamma^2 }\left(q^2- \left( 1+\xi\right)^2  \right)\chi_1(xy)+ \dfrac{1}{4\gamma^2 }\left(q^2 -\left( 1-\xi\right)^2 \right)\chi_2(xy)$$$$\\
		 		 		=\dfrac{q}{ 2\gamma   }  \chi_1(z_0) \chi_1(xy) -\dfrac{q}{ 2\gamma   } \chi_2(z_0) \chi_2(xy),
		 		 		$$
		 		 		which implies  by Lemma \ref{ui} (i) that
		 		 	$
		 		 		\dfrac{q}{ 2\gamma   }  \chi_1(z_0)=\dfrac{1}{4\gamma^2 }\left(q^2- \left( 1+\xi\right)^2  \right)
		 		 	$
		 		 		and
		 		 		$
		 		 	\dfrac{q}{ 2\gamma} \chi_2(z_0)=-\dfrac{1}{4\gamma^2 }\left(q^2 -\left( 1-\xi\right)^2  \right),$
		 		 		since $\chi_1$ and $\chi_2$ are different and non-zero. Then we deduce that
		 		 		\begin{equation*}
		 		 		\chi_1(z_0)=\dfrac{1}{2\gamma q}\left(q^2- \left( 1+\xi\right)^2  \right) \hspace{0.3cm } \text{and}\hspace{0.3cm } \chi_2(z_0)=-\dfrac{1}{2\gamma q }\left(q^2 -\left( 1-\xi\right)^2  \right).
		 		 		\end{equation*}
		 		 		So we are in  part (5).
\par (1.1.ii) There exist a  non-zero multiplicative function $\chi $ on  $S$ and an additive function $A$ on $S\setminus I_{\chi}$ such that
		 		 		$
		 		 		\gamma g= \Psi_{A\chi,\rho}$ and $\gamma f=\chi   \pm \Psi_{A\chi,\rho}.
		 		 		$
\par If $z_0 \in I_{\chi} \setminus P_{\chi}$ we have  $\gamma g(z_0)=$ $\Psi_{A\chi,\rho}(z_0)=0$ by definition of $\Psi_{A\chi,\rho}$. If $z_0 \in   P_{\chi}$ we have $\chi(z_0)=0$ and $ \lvert\gamma g(z_0)\rvert$=$\lvert\rho(z_0)\rvert$=$\lvert\chi(z_0)\pm\rho(z_0)\rvert$=$\lvert\gamma f(z_0)\rvert$$=0$. So, If $z_0 \in I_{\chi}$ we get that $\gamma g(z_0)=0$ which is a contradiction because $g(z_0)\neq0$ and $\gamma=\dfrac{g( z_0)}{f(z_0^2)}$.
		 		 		 \par Hence, $z_0 \in S \setminus I_{\chi}$ and  we have  $\chi(z_0) \neq0 $. Since  $f(z_0)=0$ by assumption  we get that
		 		 		$
		 		 		 f(z_0)=\dfrac{1}{\gamma}[  \chi(z_0)\pm A(z_0)\chi(z_0)]=0,
		 		 	 	 $
		 		 		which implies that $A(z_0)=-1  $.
		 		 		 Now   for all $x,y \in    S \setminus I_{\chi}  $, we have
		 		 		$xyz_0\in    S \setminus I_{\chi}$, then by using (\ref{4}) we get that$
		 		 		\left( \dfrac{1}{\gamma }-\chi(z_0)  \right) \chi(xy)+\left( \dfrac{1}{\gamma }+\chi(z_0)\right) \chi(xy)A(xy)=0,
		 		 		$
		 		 		which implies according to Lemma \ref{ui}(i), that
		 		 		$
		 		 		\dfrac{1}{\gamma }-\chi(z_0)=0$ and   $\dfrac{1}{\gamma }+\chi(z_0) =0,
		 		 		$
		 		 		since $A\neq0$.
		 		 		Therefore, $
		 		 		\chi(z_0)=	\dfrac{1}{\gamma }=- \dfrac{1}{\gamma },
		 		 		$
		 		 	 which is a contradiction because $	\dfrac{1}{\gamma }\neq 0$ by assumption. So we do not have a solution corresponding this this possibility.  \\
		 		 	{Subcase 1.2}:  Suppose that $ g(z_0^2) \neq0$, then the functional equation (\ref{000} ) can be rewritten as follows
		 		 		$
		 		 		\beta g(xy)=\beta^2g(x)g(y)-\beta^2f(x)f(y)+\alpha \beta f(xy), \hspace{0.3cm} x, y \in S
		 		 		$
		 		 		with
		 		 		$
		 		 		\beta:=\dfrac{g(z_0)}{g(z_0^2)}\neq0$ and $\alpha:=\dfrac{f(z_0^2)}{g(z_0^2)}.
		 		 		$
		 		 		This shows that the pair $(\beta g, \beta f)$ satisfies the functional equation (\ref{6-2}). So according to  \cite[Theorem 3.1]{h}, and taking into account that $f$ and $g$ are linearly independent, there are only the following possibilities.
		 		 		\par (1.2.i)   There exist a constant $q \in \mathbb{C} \backslash \{\pm \alpha\}$ and two different non-zero  multiplicative functions $\chi_1$ and $\chi_2$ on $S$ such that
		 		 		$
		 		 		\beta f=\alpha \dfrac{\chi_1 +\chi_2}{2}+q\dfrac{\chi_1 - \chi_2}{2} $ and $ \beta g = \dfrac{\chi_1 +\chi_2}{2}\pm \sqrt{1+q^2-\alpha^2} \dfrac{\chi_1 - \chi_2}{2}.
		 		 	$
		 		 	Introducing $\delta :=\pm \sqrt{1+q^2-\alpha^2}$ we find that
		 		 		$
		 		 	 f=\alpha \dfrac{\chi_1 +\chi_2}{2	\beta}+q\dfrac{\chi_1 - \chi_2}{2	\beta} $ and $  g = \dfrac{\chi_1 +\chi_2}{2	\beta} +\delta \dfrac{\chi_1 - \chi_2}{2	\beta}.
		 		 		$
		 		 	   By using (\ref{4}) we get
		 		 		$
		 		 		\dfrac{1}{4\beta ^2 }\left( \left( 1+\delta\right)^2 -(\alpha +q)^2\right)\chi_1(xy)+ \dfrac{1}{4\beta ^2 }\left( \left( 1-\delta\right)^2  -(\alpha +q)^2\right)\chi_2(xy)
		 		 		$	$=\dfrac{1}{ 2\beta   }\left( 1+\delta\right) \chi_1(z_0) \chi_1(xy) +\dfrac{1}{ 2\beta    }\left( 1-\delta\right) \chi_2(z_0) \chi_2(xy),$
		 		 			which implies by Lemma \ref{ui}(i) that
		 		 		$
		 		 		\dfrac{1}{ 2\beta }\left( 1+\delta \right) \chi_1(z_0)=\dfrac{1}{4\beta ^2 }  \left( (1+\delta)^2 -(\alpha +q)^2\right)
		 		 		$\\
		 		 		and\\
		 		 		$
		 		 		\dfrac{1}{ 2\beta    }\left( 1-\delta\right) \chi_2(z_0)= \dfrac{1}{4\beta ^2 }\left( \left( 1-\delta\right)^2  -(\alpha +q)^2\right),
		 		 		$
		 		 		since $\chi_1$ and $\chi_2$ are different non-zero multiplicative functions. Notice that $\delta\neq \pm 1$ because $q\neq \pm \alpha$. Therefore we deduce that
		 		 			$
		 		 		  \chi_1(z_0)=\dfrac{  ( 1+\delta )^2-(\alpha +q)^2}{2\beta(1+\delta)   }$ and  $\chi_2(z_0)=\dfrac{  (1-\delta)^2 -(\alpha +q)^2}{2\beta(1-\delta)   }.
		 		 			$
		 		 			 Hence, by writing $\gamma$ instead of $\beta$ we get part (6).
		 		 		
		 		 	 \par  (1.2.ii) $\alpha \neq0$ and there exist two different  non-zero multiplicative functions $\chi_1$ and $\chi_2$ on $S$ such that $\beta f=\alpha \chi_1$ and  $\beta g = \chi_2$. By
		 		 	    using (\ref{4}) again we get $
		 		 		\dfrac{1 }{\beta }\left(\chi_2(z_0) -\dfrac{1 }{\beta }\right)  \chi_2(xy)+\left( \dfrac{\alpha^2}{\beta^2}\right)  \chi_1(xy)=0,
		 		 	$ which gives
		 		 			$\chi_2(z_0)=\dfrac{1  }{\beta }  $     and    $  \alpha=0  $, since $\chi_1$ and $\chi_2$ are different.	So this possibility is excluded because $\alpha \neq0$.
		 		 	 \par (1.2.iii) There exist a non-zero multiplicative function $\chi $ on $S$ and an additive function $A$ on $S\setminus I_{\chi}$ such that
		 		 		$
		 		 		\beta f=\alpha \chi +\Psi_{A\chi,\rho}$ and $\beta g= \chi \pm \Psi_{A\chi,\rho},$	which gives $f=\dfrac{1}{\beta}(\alpha \chi +\Psi_{A\chi,\rho})$ and $    g=\dfrac{1}{\beta}( \chi \pm \Psi_{A\chi,\rho}).
		 		 		$
		 		 		
		 		 		\par If $g=\dfrac{1} {\beta}(\chi + \Psi_{A\chi,\rho})$ then $z_0 \notin I_{\chi} \setminus P_{\chi}$. Indeed if not we have $\chi(z_0)=0$ and $\Psi_{A\chi,\rho}(z_0)=0$. Then $\beta g(z_0)=\chi(z_0)+\Psi_{A\chi,\rho}(z_0)=0$. This   contradicts the fact that  $g(z_0)\neq0$.
		 		 		\\	On the other hand  $z_0 \notin   P_{\chi}$. Indeed if note  we have    $\chi(z_0)=0$. Then $ \beta g(z_0)=\Psi_{A\chi,\rho}(z_0)= \beta f(z_0)=0$, which is a contradiction because $g(z_0)\neq0$.
		 		 		So, $z_0 \in S \setminus I_{\chi}$ and then $\chi(z_0) \neq0 $. Since  $f(z_0)=0$ we get that
		 		 		$
		 		 		f(z_0)=\dfrac{\chi(z_0)}{\beta}[ \alpha +A(z_0)]=0,
		 		 	$
		 		 		which implies that
		 		 	$
		 		 		\label{xz}
		 		 		 A(z_0)=-\alpha.
		 		 		 $
		 		 		Now let $x,y \in S \setminus I_{\chi}$ be arbitrary. We have
		 		 		$xyz_0\in S \setminus I_{\chi}$. By using (\ref{4}) we get that
		 		 		\begin{equation}
		 		 		\label{23}
		 		 		\left( \dfrac{1-\alpha^2}{\beta^2 }+ \dfrac{\alpha-1}{\beta }\chi(z_0)\right) \chi(xy)+\left( \dfrac{1-\alpha}{\beta^2 }-\dfrac{1}{\beta }\chi(z_0)\right) \chi(xy)A(xy)=0.
		 		 		\end{equation}
		 		 			If $A=0$  then   $\rho\neq0$ because $\Psi_{A\chi,\rho}\neq0$,  $\alpha=0$, and $	\chi(z_0)=	\dfrac{1 }{\beta }$  by    (\ref{23}).  This is a special case of solution part (7).\\
		 		 			If $A\neq0$ then by Lemma \ref{ui} {(ii)}
		 		 	we get from (\ref{23}) that
		 		 		$
		 		 	 	\dfrac{1-\alpha^2}{\beta^2 }+ \dfrac{\alpha-1}{\beta }\chi(z_0) =0 $ and  $ \dfrac{1-\alpha}{\beta^2 }-\dfrac{1}{\beta }\chi(z_0)  =0.
		 		 	$
		 		 		As $\alpha\neq1$, because $\chi(z_{0})\neq0$, we deduce that
		 		 		$
		 		 		\label{xc}
		 		 		\chi(z_0)=	\dfrac{1-\alpha}{\beta }   $ and  $ 	\chi(z_0)=	\dfrac{1+\alpha}{\beta }.
		 		 		$
		 		 		So,  we obtain that $\alpha =0$ and  $\chi(z_0)=	\dfrac{1} {\beta }$,   and   the form of $f$ reduces to $f= \dfrac{1} {\beta}  \Psi_{A\chi,\rho}$.
		 		 		So we are in part  (7 ).
\par If $g=\dfrac{1} {\beta}(\chi - \Psi_{A\chi,\rho})$, by using a similar computation as above we show that we are also in part (7).\\
		 		 		{Case 2}. Suppose $f(z_0)\neq0$. By using system
		 		 		(\ref{4})  and  (\ref{03}), we deduce by an elementary computation that for any $ \lambda \in \mathbb{C}$
		 		 		\begin{equation}
		 		 		\label{80}
		 		 		(g-\lambda f)(xyz_0)=(g-\lambda f)(x )(g-\lambda f)(y)-(\lambda^2+\mu\lambda+1)f(x)f(y),\,x,y\in S.
		 		 		\end{equation}
		 		 		Let  $\lambda_1$ and $\lambda_2$  the two roots of the equation $\lambda^2+\mu\lambda+1=0$, then $\lambda_1\lambda_2=1$ which gives $\lambda_1\neq0$ and $\lambda_2\neq0$. According to \cite[Proposition 16]{i} we deduce from (\ref{80}) that
		 		 		$g-\lambda_{1} f:=\chi_1(z_0)\chi_1 \hspace{0.4cm}  and \hspace{0.4cm}  g-\lambda_{2} f:= \chi_2(z_0)\chi_2,
		 		 		$
		 		 		where $\chi_1$ and $\chi_2$ are two multiplicative functions such that $\chi_1(z_0)\neq0$  and	$\chi_2(z_0) \neq0$,    	 because $f$ and $g$ are linearly independent.
\par If $\lambda_1\neq \lambda_2 $, then   $\chi_1\neq\chi_2$  and we get that
		 		 		$
		 		 		g =\dfrac{\lambda_2\chi_1(z_0)\chi_1-\lambda_1\chi_2(z_0)\chi_2 }{\lambda_2-\lambda_1 } $  and $	f =\dfrac{ \chi_1(z_0)\chi_1- \chi_2(z_0)\chi_2 }{\lambda_2-\lambda_1 }.
		 		 		$
		 		 By putting $\lambda_1 =ic$ we get the solution categories (4).
		 		 		
\par If $\lambda_1=\lambda_2=:\lambda$, then $g-\lambda f=:\chi(z_0)\chi$ where $\chi $ is a multiplicative function on $S$ such that 	$\chi (z_0) \neq0$, because $f$ and $g$ are linearly independent. Hence
\begin{equation}\label{eq1}g=\chi(z_0)\chi+\lambda\,f.\end{equation} Substituting this (\ref{03}) an elementary computation shows that
\begin{equation*}f(xyz_0)=\chi(z_0)f(x)\chi(y) +\chi(z_0)f(y)\chi(x)+(2\lambda+\mu)f(x)f(y),\end{equation*}
for all $x,y\in S.$\\
Moreover $\lambda=1$ or $\lambda=-1$ because  $\lambda_1\lambda_2=1$. Hence, $(\lambda, \mu)=(1,-2)$ or $(\lambda, \mu)=(-1,2)$ since $\lambda^2+\mu\lambda+1=0$ and $\lambda\in\{-1,1\}$. So, the functional equation above reduces to
\begin{equation*}f(xyz_0)=\chi(z_0)f(x)\chi(y) +\chi(z_0)f(y)\chi(x),\end{equation*}
for all $x,y\in S.$ Thus the pair $(f, \chi)$  satisfies  the functional equation (\ref{bb}). Hence, in  view of Proposition \ref{t2}, we get  $f=A(z_0)\chi +\Psi_{A\chi,\rho}.$ Then, by using (\ref{eq1}), we derive that $g=\chi (z_0) \chi +\lambda\,f= (\chi(z_0)+\lambda\,A(z_0))\chi+\lambda^{2}\Psi_{A\chi,\rho}=(\chi(z_0)\pm A(z_0))\chi+\Psi_{A\chi,\rho}$. This is part (8).
\par Conversely it is easy to check that the formulas of $f$ and $g$ listed in  Theorem  \ref{jh} define solutions of (\ref{4}).
\par Finally, suppose that  $S$ is a topological semigroup. The continuity of the solutions of the forms (1)-(6) follows directly from \cite[Theorem 3.18]{g}, and for the ones of the forms (7) and (8)  it is parallel  to the proof used in \cite[Theorem 2.1]{b} for categories (7)-(8). This completes the proof of Theorem \ref{jh}.
\end{proof}
           	
\end{document}